\input amstex
\TagsOnRight
\documentstyle{amsppt}
\input epsf
\def\epsfsize#1#2{\hsize}
\def\fig#1{\smallskip\centerline{\epsffile{FIG#1.EPS}}\centerline{\bf Fig\. #1}\bigskip}

\let\but\setminus \let\emb\hookrightarrow 
\let\x\times \let\tl\tilde \let\eps\varepsilon
\def\G{\Cal G} \def\P{\Cal P} \def\M{\Cal M} \def\R{\Bbb R} \def\Z{\Bbb Z}

\topmatter 

\title $n$-quasi-isotopy: III. Engel conditions \endtitle
\author Sergey A. Melikhov and Roman V. Mikhailov \endauthor

\subjclass Primary: 57M25; Secondary: 57M05, 57Q30,
57Q37, 20F12, 20F14, 20F34, 20F45 \endsubjclass
\keywords $k$-quasi-isotopy, link homotopy, PL isotopy, Milnor link
group,
Milnor $\bar\mu$-invariants, generalized Sato--Levine invariant, Cochran
derived invariants, Casson--Walker--Lescop invariant, basic commutators,
bounded left Engel elements, subnormal cyclic subgroups
\endkeywords
\thanks Partially supported by the Russian Foundation for Basic Research
Grants 99-01-00009 and 02-01-00014 and the Latvian Science Society Grant
No.\ 01.0282 \endthanks
\address Steklov Mathematical Institute, Division of Geometry and Topology;
\newline $\Gamma$C$\Pi$-1, ul. Gubkina 8, Moscow 117966, Russia \endaddress
\email sergey\@melikhov.mccme.ru, rmikhailov\@mail.ru \endemail
\address{\it First author's current address:}
University of Florida, Department of Mathematics;\newline
358 Little Hall, PO Box 118105, Gainesville, FL 32611-8105, U.S. \endaddress
\email melikhov\@math.ufl.edu \endemail

\abstract In part I it was shown that for each $k\ge 1$ the generalized
Sato--Levine invariant detects a gap between $k$-quasi-isotopy of link
and peripheral structure preserving isomorphism of the finest quotient $\G_k$
of its fundamental group, `functorially' invariant under $k$-quasi-isotopy.
Here we show that Cochran's derived invariant $\beta^k$, provided $k\ge 3$, and
a series of $\bar\mu$-invariants, starting with $\bar\mu(111112122)$ for $k=3$,
also fall in this gap.
In fact, all $\bar\mu$-invariants where each index occurs at most $k+1$ times,
except perhaps for one occuring $k+2$ times, can be extracted from $\G_k$, and
if they vanish, $\G_k$ is the same as that of the unlink.

We also study the equivalence relation on links (called `fine
$k$-quasi-isotopy') generated by ambient isotopy and the operation of interior
connected sum with the second component of the $(k+1)^{\text{th}}$ Milnor's
link, where the complement to its first component is embedded into the link
complement.
We show that the finest quotient of the fundamental group, functorially
invariant under fine $k$-quasi-isotopy, is obtained from the fundamental group
by forcing all meridians to be $(k+2)$-Engel elements.
We prove that any group generated by two $3$-Engel elements has lower central
series of length $\le 5$.
\endabstract
\endtopmatter

\document 

\head 1. Introduction \endhead

Two links $L_1,L_2\:mS^1\emb S^3$, where $mS^1=S^1_1\sqcup\dots\sqcup S^1_m$,
are called {\it $k$-quasi-isotopic} if they are PL homotopic in the class of
{\it $k$-quasi-embeddings} with at most one double point, defined as follows.
Any PL embedding is set to be a $k$-quasi-embedding, for arbitrary $k$.
A PL map $f\:mS^1\to S^3$ with a single double point
$f(p)=f(q)$ is said to be a $k$-quasi-embedding if $p,q$ lie in the same
component $S^1_i$ and (in the case $k>0$) moreover, in addition to
the singleton $P_0=\{f(p)\}$, there exist subpolyhedra
$P_1\i\dots\i P_k\i S^3\but\bigcup_{j\ne i}f(S^1_j)$ and arcs
$J_0\i\dots\i J_{k-1}\i S^1_i$ such that $f(J_j)\i P_{j+1}$ and
$f^{-1}(P_j)\i J_j$, moreover the inclusion $P_j\cup f(J_j)\i P_{j+1}$
is null-homotopic in $P_{j+1}$, for each $j=0,\dots,k-1$.
(Beware that in \cite{MM} $k$-quasi-isotopy was what we call strong
$k$-quasi-isotopy now.)
The {\it $k$-inessential lobe} of the $k$-quasi-embedding $f$ is
$f(\tilde J_0)$, where $\tilde J_0$ is the subarc of $J_0$ such that
$\partial\tilde J_0=\{p,q\}$.

This definition was suggested by certain well-known higher dimensional
constructions, in particular the Penrose--Whitehead--Zeeman--Irwin trick and
the Casson handles, see \cite{MR1} for details and references.
The importance of these equivalence relations lies in that they all resemble
link homotopy in satisfying the following property: given a link in $S^3$,
possibly wild, there exists an $\eps=\eps(k)>0$ such that all $\eps$-close PL
links are $k$-quasi-isotopic to each other \cite{MR1}.

It was observed in \cite{MR1} that the $k^{\text{th}}$ Milnor's link
\cite{Mi2; Fig. ~1} (see Fig. ~5 below, where $k=4$) is $(k-1)$-quasi-isotopic
to the unlink but not $k$-quasi-isotopic to it.
Also the $k^{\text{th}}$ Whitehead link $\Cal W_k$ (i.e\. the $k$-fold
untwisted left-handed Whitehead double of the Hopf link, where the doubling
is performed on either one of the components at each stage) is
$(k-1)$-quasi-isotopic to the unlink (see Figures 4 and 5 in \cite{MR1}).
We do not presently know any invariants that could tell us whether the boundary
link $\Cal W_k$ is not $k$-quasi-isotopically trivial for $k\ge 2$.
Also it is still not clear to us whether connected sum `accumulates complexity'
up to $k$-quasi-isotopy, $k\ge 1$, in the sense of Problem 0.2 from \cite{MR1},
which was the initial motivation of the present series of papers.

It was proved in \cite{MR1}, \cite{MR2} that the following are invariant under
$k$-quasi-isotopy: Cochran's derived invariants $\beta^i$, $i\le k$; finite
type invariants of type $\le k$ (either in the usual sense, or more generally
in the sense of Kirk--Livingston \cite{KL}) that are well-defined up to PL
isotopy; and Milnor's $\bar\mu$-invariants of length $\le 2k+3$.
It follows that any linear combination of Conway's polynomial coefficients at
the powers $\le k+m-1$, invariant under PL isotopy, as well as the first
non-vanishing coefficient of the Conway polynomial, if it occurs at a
power less than $(2k+3)(m-1)$, are invariant under $k$-quasi-isotopy.
In this paper we also show that all $\bar\mu$-invariants with no more than
$k+1$ entries of each index are invariant under $k$-quasi-isotopy.

The referee suggested us to mention a connection with capped half-gropes, and
conjectured a version of the part (a) of the following proposition.
Let $M$ be an oriented $2$-manifold; if $M$ is a disk, it will be omitted
from the notation.
An {\it $M$-like capped half-grope of class $1$} is a polyhedral pair
$(G_1,\partial G_1)$ obtained from $(M,\partial M)$ by choosing an even number
of disjoint disks $D^2_1,\dots,D^2_{2n}$ in the interior of $M$ and identifying
each $D^2_{2i-1}$ with $D^2_{2i}$ by an orientation-preserving homeomorphism
$h_i$.
The images of $D_i$ under the identification $M\to G_1$ are the {\it caps},
the complement in $G_1$ to their interiors is the {\it surface}, and the images
in $G_1$ of any collection of $n$ disjoint paths in
$M\but (D_1\cup\dots\cup D_{2n})$, connecting a point $p_i\in\partial D_{2i-1}$
to $h_i(p_i)\in\partial D_{2i}$ for each $i$, will be called an {\it expansion
base} for $G_1$.
An {\it $M$-like capped half-grope of class $k+1$} is a polyhedral pair
$(G_{k+1},\partial G_{k+1})$ obtained from $(G_1,\partial G_1)$ by attaching
a (disk-like) capped half-grope $G_k(i)$ of class $k$ along its boundary
$\partial G_k(i)$ to each circle $S^1_i$ in some expansion base for $G_1$.
{\it Caps} and {\it surfaces} of $G_{k+1}$ are those of $G_1$ and of each
$G_k(i)$.

\proclaim{Proposition 1.1}
a) Two $m$-component links related by $k$-quasi-isotopy cobound a link map of
$m$ annulus-like capped half-gropes of class $k+1$ into $S^3\x I$, which embeds
the bottom stage surface and the union of all other surfaces of each grope.

b) The $k$-inessential lobe of a $k$-quasi-embedding $f\:mS^1\to S^3$ with
a single double point on the $i^{\text{th}}$ component bounds an immersion of
a capped half-grope of class $k$ into $S^3\but\bigcup_{j\ne i}f(S^1_j)$ which
embeds each surface of the grope, moreover any two surfaces from the same stage
meet tranversely in their interiors.
\endproclaim

The converse to (a) is false, since concordance does not imply
$1$-quasi-isotopy \cite{MR1}.

\demo{Proof}
To prove (b), convert the track of the null-homotopy of the $k$-inessential
lobe in $P_1$ into an embedded surface (by tricks II and IV from Appendix A of
Cochran's memoire \cite{Co}), then (if $k>1$) choose an expansion base so that
any two circles have zero linking number (this can be achieved by dragging one
end of a handle between the ends of another handle), and use the same tricks to
convert the tracks of their transverse null-homotopies in $P_2$ into embedded
surfaces, meeting only in their interiors.
This process can be continued; in order to choose expansion bases for different
class $1$ capped half-gropes from the same stage so that circles from one have
trivial linking numbers with that from the other, it may be necessary to go
back to lower stages and use handle slides.
To prove (a), convert the immersed annuli, forming the image of
the $k$-quasi-isotopy in $S^3\x I$, into surfaces by resolving
self-intersections, and then perturb the surfaces of each grope from (b) in
the fourth coordinate, so as to place them consecutively in an order
of increasing stage number. \qed
\enddemo

We define $\G_k(L)$ to be the quotient of the fundamental group
$\pi(L)$ by the normal subgroup
$$\mu_k=\left<[m,m^g]\mid m\in M,g\in\underbrace{\left<m\right>^{.^{.^{.^{
\left<\!m\!\right>^{\pi(L)}}}}}\hskip-18pt}_{\text{$k$ times}}\hskip18pt\right>,$$
where $\left<m\right>^H$ denotes the normal closure of $m$ in the subgroup $H$,
and $M$ denotes the set of all meridians.
\footnote{We use the left normed notation $h^g=g^{-1}hg$,
$[g,h]=g^{-1}h^{-1}gh$ throughout the paper.}
Notice that for $k=0$ this coincides with Milnor's link group $\G(L)$
\cite{Mi1} which is not surprising as $0$-quasi-isotopy is just link homotopy.
Equivalently, $\mu_k$ can be described as the normal closure of
$$\bigcup_{i=1}^m\left<[m_i,m_i^g]\mid g\in\left<m_i\right>_k^{\pi(L)}\right>,$$
where each $m_i$ is a fixed meridian to the $i^{\text{th}}$ component and
$H_k^G$ denotes $H^{H_{k-1}^G}$ where $H_0^G=G$.
In \S2 we will give more descriptions of $\mu_k$.
It was shown in \cite{MR1, Theorems 3.2 and 3.7} that $\G_k$ is the finest
quotient of $\pi(L)$, invariant under $k$-quasi-isotopy {\it functorially},
i.e\. so that for any links $L_0$, $L_1$ obtained one from another by an
allowed crossing change with $L_s$ being the intermediate singular link,
the diagram
$$\CD
\pi_1(S^3\but\text{regular neigh}\hskip -20pt@.\hskip-20pt
\text{borhood of }L_s(mS^1))\\
@Vi_*VV@VVi_*V\\
\pi(L_0)@.\pi(L_1)\\
@Vp_{L_0}VV@VVp_{L_1}V\\
\G_k(L_0)@>\simeq>>\G_k(L_1)
\endCD$$
commutes.
Let $\P_k$ be the finest peripheral structure in $\G_k(L)$, invariant under
$k$-quasi-isotopy, more precisely, $\P_k$ is the collection of $m$ pairs
$(\bar m_i,\Lambda_i)$, each defined up to simultaneous conjugation, where
$\bar m_i\in\G_k(L)$ is the coset of a meridian $m_i\in\G_k(L)$ to the
$i^{\text{th}}$ component of $L$, and $\Lambda_i\i\G_k(L)$ is the set of
cosets $\bar l_{i\alpha}\in\G_k(L_\alpha)\simeq\G_k(L)$ of the longitudes
$l_{i\alpha}$ corresponding to some representatives $m_{i\alpha}$ of $\bar m_i$
in the fundamental groups of all links $L_\alpha$, $k$-quasi-isotopic to $L$
(that $m_{i\alpha}$ are meridians follows from functoriality).
It follows from \cite{MR1; proofs of Theorems 3.2 and 3.7} that the set
$\Lambda_i$ is the right coset $N_k(m_i)\bar l_i$ of the subgroup
$$N_k(m_i)=\left<[g^{-1},g^m]\mid g\in\left<m\right>_k^{\G_k}\right>$$
containing the class $\bar l_i$ of the longitude $l_i$ corresponding to $m_i$.

\proclaim{Theorem 1.2}
a) \cite{MR1} If $k\ge 1$, the generalized Sato--Levine invariant is invariant
under $k$-quasi-isotopy but cannot be extracted from $(\G_k,\P_k)$.

b) If $k\ge 3$, Cochran's derived invariants
$\beta^{[k/2]+2}_\pm,\dots,\beta^k_\pm$ (see \cite{Co}) are invariant under
$k$-quasi-isotopy but cannot be extracted from $(\G_k,\P_k)$.

c) If $k\ge 3$, a nonzero $\bar\mu$-invariant of $2$-component links of length
$\le 2k+3$ with at least $k+3$ entries of one of the indices is invariant
under $k$-quasi-isotopy but cannot be extracted from $(\G_k,\P_k)$.
\endproclaim

\remark{Remark} Of course, $k\ge 3$ is no real restriction in (c), since there
are no nontrivial $\bar\mu$-invariants of two-component links of length $5$ or
$7$, and those of length $4$ or $6$ must contain at least two entries of each
index.
Notice also that the only invariants of length $8$ that satisfy the conditions
of (c) with $k=3$ are $\bar\mu(11111122)$ and $\bar\mu(11222222)$ \cite{Mi2},
which in the case of zero linking number coincide with the residue classes of
Cochran's $\beta^3_{\pm}$ mod g.c.d\. of the Sato-Levine invariant and
$\beta^2_{\pm}$ \cite{Co}.
But there are $4$ linearly independent $\bar\mu$-invariants of $2$-component
links of length $9$ \cite{Orr}, and a computation using Proposition 2.13 below
along with the second Witt formula (for the number of basic commutators with
prescribed entries) \cite{Hall; eq\. (11.2.4)}, \cite{MKS} shows that one of
these has $6$ entries of the index $``1"$, thus satisfying the hypothesis of
(c) with $k=3$.
Using Milnor's symmetry relations, it is easy to verify that this invariant is
$\frac15\bar\mu(112111122)=\frac12\bar\mu(111112122)$.
\endremark

\demo{Proof of 1.2}
The part (c) follows from Theorem 2.12c and \cite{MR1; Corollary 3.5}.
The invariance in part (b) was proved in \cite{MR2}, and since $\beta^k$ is
an integer lifting of one of the invariants from (c), it cannot be extracted.
\qed
\enddemo

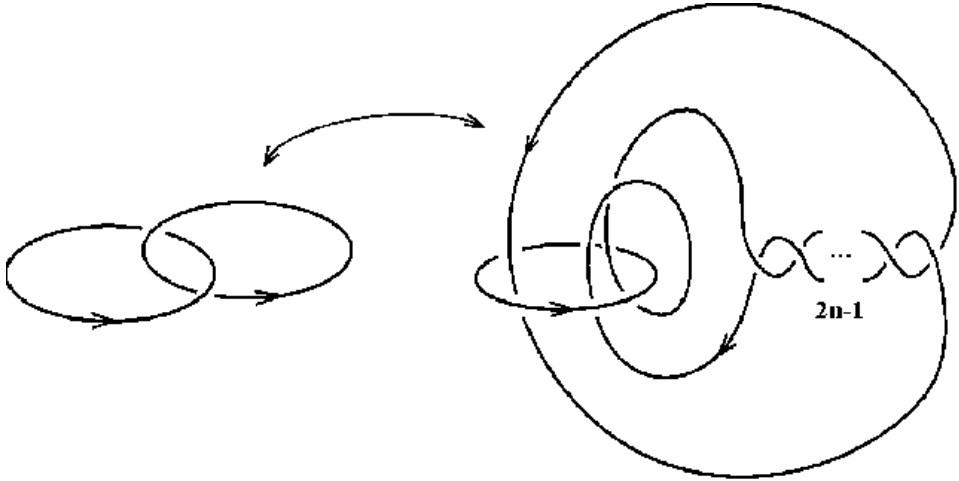
\fig 1

The following argument is a modification of \cite{MR1; proof of Corollary 1.3}.

\example{Example 1.3} The infinite family of links $\M_n$, depicted on the
right-hand side of Fig\. 1, has the same $(\G_k,\P_k)$ as the Hopf link (for
each $k$), yet no two links of this family are $1$-quasi-isotopic.

To see the latter, let us recall the axiomatic definition of the generalized
Sato--Levine invariant $\tl\beta$, constructed independently by Polyak--Viro,
Kirk--Livingston, and Akhmet'ev (see \cite{KL}, \cite{AR}, \cite{AMR}).

Given two links $L_0,L_1\:S^1\sqcup S^1\emb\R^3$, related by a single crossing
change in one of the components, of the type shown on Fig\. 2a, so that the
two lobes of the singular component of the intermediate singular link have
linking numbers $n$ and $l-n$ with the other component, where $l$ is the
linking number of $L_0$, one postulates the formula
$\tl\beta(L_0)-\tl\beta(L_1)=n(l-n)$.
In addition, it is assumed that $\tl\beta(\Cal H_n)=0$, where $\Cal H_n$ is
shown on Fig\. 2b.
It was proved in \cite{AR}, \cite{KL} that these two axioms yield
a well-defined invariant of ambient isotopy.
Moreover, it is easy to see that this invariant $\tl\beta$ is unchanged under
$1$-quasi-isotopy.

\fig 2

We note that $\M_n$ can be obtained from $\M_{n-1}$ by one crossing
change, where the jump of $\tl\beta$ is $-2$.
Hence $\tl\beta(\M_n)=-2n$ which proves that these links are pairwise
non-$1$-quasi-isotopic.

On the other hand, it is not hard to see from the Wirtinger presentation that
each fundamental group $\pi(\M_n)$ is generated by two meridians, say
$m_1$ and $m_2$.
Then $[m_1,l_1]=1$ in $\pi=\pi(\M_n)$, where $l_1$ is the longitude
corresponding to $m_1$.
Since the linking number is $1$, we have $l_1=cm_2$ for some $c\in[\pi,\pi]$.
Hence $1=[m_1,cm_2]=[m_1,m_2][m_1,c]^{m_2}$ and so
$[m_1,m_2]=[c,m_1]^{m_2}\in [[\pi,\pi],\pi]$.
Thus $\gamma_2\pi=\gamma_3\pi$ and it follows that $\pi$ has length $2$, i.e\.
$\gamma_k\pi=\gamma_2\pi$ for any $k\ge 2$.
(We recall that the finite lower central series of $\pi$ is defined inductively
by $\gamma_1\pi=\pi$ and $\gamma_{n+1}\pi=[\gamma_n\pi,\pi]$.)
But for every link $L$, the group $\G_k(L)$ is nilpotent by \cite{MR1}
(or by Theorem 2.1 below).
Therefore $\G_k(\M_n)$ is Abelian, more precisely it is $\Z\oplus\Z$.
In particular, $\P_k$ coincides with the image of the peripheral structure
in $\pi$ under the quotient map $\pi\to\G_k$. \qed
\endexample

\example{Example 1.4} Let us point out certain multi-component analogues of
the generalized Sato--Levine invariant, well-defined up to $1$-quasi-isotopy.

For each rational vector $s=(s_1,\dots,s_m)\in\Bbb Q^m$ and any choice of
a base link in each link homotopy class of $m$-component links, the following
crossing change formula yields a well-defined link invariant $\beta(L,s)$
\cite{J2}:
$$\beta(L_1,s)-\beta(L_0,s)=\pmatrix
s_1       &\hdots &a_{1,x-1}   &q_1       &a_{1,x+1}   &\hdots &a_{1m}    \\
\vdots    &\ddots &\vdots      &\vdots    &\vdots      &\ddots &\vdots    \\
a_{x-1,1} &\hdots &s_{x-1}     &q_{x-1}   &a_{x-1,x+1} &\hdots &a_{x-1,m} \\
p_1       &\hdots &p_{x-1}     &\lambda   &p_{x+1}     &\hdots &p_m       \\
a_{x+1,1} &\hdots &a_{x+1,x-1} &q_{x+1}   &s_{x+1}     &\hdots &a_{x+1,m} \\
\vdots    &\ddots &\vdots      &\vdots    &\vdots      &\ddots &\vdots    \\
a_{m1}    &\hdots &a_{m,x-1}   &q_m       &a_{m,x+1}   &\hdots &s_m
\endpmatrix$$
Here $L_0$ and $L_1$ differ by a crossing change in the $x^{\text{th}}$
component, of the type shown on Fig\. 2a, $a_{ij}$ are the linking numbers
between the components, $p_i$ (resp\. $q_i$) are the linking numbers between
the first (resp\. second) lobe of the singular $x^{\text{th}}$ component and
the other components (so that $p_i+q_i=a_{xi}$), and $\lambda$ is the linking
number between the two smoothed lobes.

(In fact, let $M_s$ denote the $3$-manifold obtained from $S^3$ by the surgery
on $(L,s)$, and $A=(a_{ij})$ denote the surgery matrix, where $a_{ii}=s_i$.
Then up to a constant depending on the choice of base links, $\beta(L,s)$ is
the Casson--Walker--Lescop invariant of $M_s$ divided by $(-1)^bq_1\dots q_m$
where $b$ is the number of negative eigenvalues of $A$ and
$s_i=\frac{p_i}{q_i}$ with $(p_i,q_i)=1$, $q_i>0$ \cite{J2}.
Also if $\det A\ne 0$, then up to the same constant, $\beta(L,s)$ is the
Casson--Walker invariant of the rational homology sphere $M_s$ divided by
$\frac2{\det A}$ \cite{J1}.)

We note that $\beta(L,s)$ is invariant under $1$-quasi-isotopy if $s$ is such
that for each $i$ the minor $\det A_i=0$, where $A_i$ is obtained from $A_i$ by
deleting the $i^{\text{th}}$ line and the $i^{\text{th}}$ column.
If $m=2$ this means that $s=(0,0)$, thus the only invariant of
$1$-quasi-isotopy yielded by this construction for $2$-component links is the
generalized Sato--Levine invariant (up to a constant).
If $m=3$ and all $a_{ij}=0$, again it is not hard to see that if $\beta(L,s)$
is invariant under $1$-quasi-isotopy, it is a linear function with constant
coefficients of the generalized Sato--Levine invariant of a two-component
sublink.

\fig 3

Let us consider the case $m=3$, all $a_{ij}\ne 0$ in more detail.
Clearly, there are only two ways of choosing $s_i$'s so that all
minors $\det A_i=0$, namely, $s_i=\pm\frac{a_{ij}a_{ik}}{a_{jk}}$.
The positive sign implies $\det A=0$, while the negative sign
leads to $\det A\ne 0$. Thus we obtain two invariants $\beta_+$,
$\beta_-$ of $1$-quasi-isotopy of $3$-component links whose
two-component sublinks have nonzero linking number. A
straightforward computation shows that each of these two
invariants distinguishes two links, shown on Fig\. 3, which have
their $2$-component sublinks ambient isotopic.
\endexample

\example{Example 1.5} It is not difficult to construct
$2$-component links with $(\G_1,\P_1)$ being the same as that of
the unlink, for each $k$, but not obviously $1$-quasi-isotopic to
the unlink; see for example \cite{MR1, Remark in the end of \S1}.
Here is a particular family of such candidates (one can verify
that all $\bar\mu$-invariants of these links of length $\le 5$
vanish, using the skein relation for the Sato--Levine invariant
$\bar\mu(1122)$).
\endexample

\fig 4

We acknowledge P. ~M. ~Akhmet'ev, K. ~K. ~Andreev, T. ~Cochran,
A. ~N. ~Dranishnikov, A. ~Ju\. ~Ol'shanskij, B. ~I. ~Plotkin,
Ju\. ~P. ~Razmyslov, R. ~Sadykov, V. ~Shpilrain, A. ~Turull for stimulating
discussions and the referee for useful remarks.
A preliminary report on the present research can be found in \cite{MM}.
Theorems 3.9 and 3.7, Proposition 3.1 and Examples 1.3--1.5 are due to
the second author; the results of \S2, Theorems 3.2 and 3.6 are due to
the first author.

\head 2. The structure of $\G_k(L)$ \endhead

Let us fix some notation.
Unless otherwise mentioned, we work in the PL category.
$F_m$ will denote the free group on $m$ fixed generators, identified with a
fixed set of meridians of the trivial $m$-component link.
By $\gamma_n G$ (or just $\gamma_n$ if $G$ is clear from the context) we denote
the $n^{\text{th}}$ term of the finite lower central series of the group $G$,
defined inductively by $\gamma_1=G$ and $\gamma_{n+1}=[\gamma_n,G]$.
We use the left-normed convention for multiple commutators:
$[a_1,\dots,a_n]=[[\dots[a_1,a_2],\dots,a_{n-1}],a_n]$.
We also use the notation $g^{-h}=(g^h)^{-1}$.

It was proved in \cite{MR1} that $\G_k(L)$ is nilpotent for any $k$ and $L$.
The same approach, based on the Hirsch--Plotkin Theorem, can also be used to
obtain an upper bound estimate of the nilpotent class of $\G_k(L)$
(depending on the minimal number of meridianal generators of $\pi(L)$).
However, this estimate turns out to be rather weak, for example, using this
method, even with some improvements, we could obtain only class $\le 6$ for
$\G_1(L)$ where $\pi(L)$ is generated by two meridians.

\proclaim{Theorem 2.1} If $L$ has $m$ components, $\G_k(L)$ is nilpotent of
class $m(k+1)$.
\endproclaim

This is immediate from the following three easy lemmas.

\proclaim{Lemma 2.2} \cite{Chen}, \cite{MKS; Lemma 5.9} Given a choice of
meridians, one to each component of $L$, then every nilpotent quotient of
$\pi(L)$ is generated by their images.
\endproclaim

(Indeed, any generator of the Wirtinger presentation is conjugate to one of the
fixed meridians, so any element of $\gamma_i$, which modulo $\gamma_{i+1}$ is a
product of commutators of weight $i$ in the generators, can be expressed modulo
$\gamma_{i+1}$ in the commutators of weight $i$ in the fixed meridians; this
follows using the commutator identities $[ab,c]=[a,c]^b[b,c]$,
$[a,bc]=[a,c][a,b]^c$ and $x^y=x[x,y]$.)

\proclaim{Lemma 2.3} Let $G$ be a group generated by $x_1,\dots,x_m$.
Then $\gamma_n G$ is the normal closure of all left-normed commutators
$[x_{i_1},\dots,x_{i_n}]$.
\endproclaim

The proof is similar to the proofs of \cite{Hall; Theorem 10.2.3},
\cite{MKS; Theorem ~5.4} and is given below for completeness.

\proclaim{Lemma 2.4} In the group $\G_k(L)$, any left-normed commutator
$[g_1,\dots,g_n]$ with at least $k+2$ entries $g_{i_1},\dots,g_{i_{k+2}}$
occupied by certain meridian\footnote{The images of meridians in
$\G_k(L)$ will also be called meridians.}, is trivial.
\endproclaim

\demo{Proof} Notice that for any elements $m,g$ of a group $G$, the commutator
$[m,g]=m^{-1}m^g$ lies in the normal closure $\left<m\right>^G$, and the same
holds for $[g,m]$.
Moreover, even if $m$ is replaced by some $m'\in\left<m\right>^G$, we still
have $[m',g]\in\left<m\right>^G$.

Consequently, if the meridian in question, which will also be denoted by $m$,
first occurs in position $i_1$, the commutator
$[g_1,\dots,g_{i_1}]=[[g_1,\dots,g_{i_1-1}],m]$ will be in
$\left<m\right>^{\G_k}$, and moreover all the subsequent commutators
$[g_1,\dots,g_{i_1+j}]$, where $j\ge 1$, will also lie
$\left<m\right>^{\G_k}$.
Now if $m$ occurs also in position $i_2$, the same reasoning ensures that the
commutator $[g_1,\dots,g_{i_2}]$ must be in
$\left<m\right>^{\left<m\right>^{\G_k}}=\left<m\right>_2^{\G_k}$,
and in the same fashion it follows that
$[g_1,\dots,g_{i_{k+2}-1}]\in\left<m\right>_{k+1}^{\G_k}$.
But $m$ is in the center of $\left<m\right>_{k+1}^{\G_k}$, according to
\cite{MR1; formula (3.1)}.
Therefore the next bracket $[g_1,\dots,g_{i_{k+2}-1},m]$ is trivial, and so
the whole commutator must be trivial. \qed
\enddemo

\demo{Proof of 2.3}
By the definition, any element of $\gamma_n G$ can be presented as
$\prod_{i=1}^r[d_i,g_i]$ where $d_i\in\gamma_{n-1} G$ and $g_i\in G$.
By induction, we can assume that each $d_i$ is a product of conjugates of
left-normed commutators of weight $n-1$ in $x_i$'s and of their inverses,
so that $[d_i,g_i]=[c_1^{\pm h_1}\dots c_s^{\pm h_s},g_i]=[c_1^{\pm h_1},v_i]
^{c_2^{\pm h_2}\dots c_s^{\pm h_s}}\dots[c_s^{\pm h_s},v_i]$,
where one uses the commutator identity $[ab,c]=[a,c]^b[b,c]$.
Now for each $j$ we have
$[c_j^{\pm h_j},v_i]=[c_j^{\pm 1},v_i^{h_j^{-1}}]^{h_j}$,
and if the sign is negative, we can eliminate it using the identity
$[a^{-1},b]=[a,b]^{-a^{-1}}$.
Analogously, $[c_j,v_j^{h_j^{-1}}]=[c_j,\prod_{k=1}^t x_{i_k}^{\pm 1}]$ can be
written as a product of conjugates of the commutators $[c_j,x_{i_k}]^{\pm 1}$.
\qed
\enddemo

\proclaim{Proposition 2.5} Any multiple commutator in the group $\G_k(L)$
with at least $k+2$ entries occupied by certain meridian, is trivial.
\endproclaim

\demo{Proof} To each multiple commutator there corresponds a finite directed
binary tree with vertices representing the bracket pairs and edges going
in the direction of (the vertex representing) the next exterior bracket.
Let $L$ denote a directed path of maximal length in this graph, and let us
define {\it complexity} of the commutator as the number of vertices not
belonging to this path.

A commutator of complexity zero is one where no two brackets are separated by
a comma, or in other words where the pairs of brackets form a chain.
The proof of Lemma 2.4 works as well in this case (alternatively, one could
reduce the case of complexity zero to the left-normed case similarly to what
follows).

Let us assume that the assertion holds for all commutators of complexity $<c$,
and let us consider a commutator of complexity $c>0$ and weight $n>3$.
By downward induction on the weight (starting from the nilpotent class of
$\G_k$) we may assume also that the statement of the theorem is true
for all commutators of weight $>n$.
Since $c>0$, there is a bracket not in the chain $L$, moreover we can assume
that the next exterior bracket is in $L$.
This must look like $$\dots [[a,b],c]\dots\qquad\text{ or }
\qquad\dots [c,[b,a]]\dots$$ where $a$ and $b$ denote either commutators or
single entries, and $c$ necessarily denotes a commutator since the exterior
bracket lies in the chain $L$.
In fact these two configurations are equivalent modulo commutators of higher
weight but still containing at least $k+2$ occurences of our meridian,
which follows using the commutator identities $[x^{-1},y]=[x,y]^{-x}$,
$x^y=x[x,y]$ and $[ab,c]=[a,c]^b[b,c]$.
Now in the situation of the first configuration, apply the Hall--Witt identity
$$[[a,b^{-1}],c]^b[[c,a^{-1}],b]^a[[b,c^{-1}],a]^c=1$$ which modulo commutators
of higher weight can be simplified as $$[[a,b],c]\sim [[a,c],b][[c,b],a].$$
Clearly, in our case these commutators of higher weight are going to retain
the $\ge k+2$ occurences of the meridian, so our commutator of weight $n$ can
be expressed as the product of two new commutators of weight $n$ where the
expression $[[a,b],c]$ is substituted respectively by $[[a,c],b]$ and
$[[c,b],a]$.
Notice that in both cases the commutator $c$ has moved into a deeper bracket,
so that the length of the chain $L$ has increased by one (in each of the two
trees corresponding to the new commutators).
It remains to apply the inductive hypothesis. \qed
\enddemo

\proclaim{Lemma 2.6} \cite{Chan; Lemma 2} Let $G$ be a group and $g$ an
element of $G$.
Then $[g,g_n^G]=[G,_{n+1}\! g]\underset\text{def}\to=
\underbrace{[[\dots[}_{n+1} G,g]\dots,g],g]$.
\endproclaim

As we could not find an English translation of \cite{Chan}, we provide a short
proof for convenience of the reader.

\demo{Proof} By induction on $n$.
It is immediate from the definition that $$\underbrace{[[\dots[}_{n+1}
G,g]\dots,g],g]=[g,[g,\dots[g,G\underbrace{]\dots]]}_{n+1}.$$
Therefore it remains to prove that $[g,[g,g_{n-1}^G]]=[g,g_n^G]$.
Indeed,
$$[g,[g,h_1]^{\alpha_1}\dots [g,h_r]^{\alpha_r}]=[g,(g^{-1}g^{h_1})^{\alpha_1}
\dots (g^{-1}g^{h_r})^{\alpha_r}].$$
For some $m_0,\dots,m_r\in\{-1,0,1\}$, the latter expression can be rewritten
as $$\multline
[g,g^{m_0}(g^{h_1})^{\alpha_1}\dots g^{m_{r-1}}(g^{h_r})^{\alpha_r}g^{m_r}]=
[g,g^{m_0}(g^{h_1})^{\alpha_1}\dots g^{m_{r-1}+m_r}(g^{h_r g^{m_r}})
^{\alpha_r}]\\
=\dots=[g,g^{m_0+\dots+m_r}(g^{h_1 g^{m_1+\dots+m_r}})^{\alpha_1}\dots
(g^{h_r g^{m_r}})^{\alpha_r}]=
[g,(g^{h'_1})^{\alpha_1}\dots (g^{h'_r})^{\alpha_r}]
\endmultline$$
for some new $h'_1,\dots,h_r'\in g_{n-1}^G$. \qed
\enddemo

Milnor's original definition of $\G(L)$ was based on the trivial case $k=0$
of the following fact (compare \cite{MR1; remark in the end of \S3}).

\proclaim{Theorem 2.7} $\mu_k=(\gamma_{k+2} A_1)\dots(\gamma_{k+2} A_m)$
where $A_i$ denotes the normal closure of a meridian to the $i^{\text{th}}$
component.
\endproclaim

\demo{Proof} The inclusion `$\i$' follows from \cite{MR1; formula (3.1)} and
Lemma 2.6.
To prove the reverse inclusion, notice that any $(k+2)$-fold commutator in
elements of $A_i$ can be rewritten, using the commutator identities, as
a product of commutators where a fixed meridian $m_i\in A_i$ occurs at least
$k+2$ times.
By Proposition 2.5, these commutators lie in $\mu_k$. \qed
\enddemo

For $k=0$ the following was proved by Levine \cite{Le}:

\proclaim{Theorem 2.8} For $i=1,\dots,m$, fix a meridian $m_i$ to the
$i^{\text{th}}$ component of $L$.
Then $\mu_k$ is the subgroup generated by $\gamma_{m(k+1)+1}$ and the basic
commutators of weight $\le m(k+1)$ in $m_i$'s with at least $k+2$ occurences
of some $m_j$.
\endproclaim

We refer to \cite{Hall} for a concise treatment of basic commutators
(see also \cite{MKS}).

\demo{Proof} A half of the statement is contained in Theorem 2.1 and
Proposition 2.5.
To prove the other half, by \cite{MR1; formula (3.1)}, Lemma 2.6 and a
straightforward induction on weight (using that $\gamma_r/\gamma_{r+2}$ is
Abelian for $r>1$) it suffices to show that for each $r\le m(k+1)$, any element
$c$ of the subgroup $[\pi(L),_{k+2}\left<m_0\right>]\cap\gamma_r\pi(L)$,
where $m_0$ is a meridian, is expressible modulo $\gamma_{r+1}\pi(L)$ in basic
commutators of weight $r$ in $m_i$'s with at least $k+2$ occurences of some
$m_j$.

We may assume without loss of generality that $\pi(L)$ is the free group $F_m$.
If $m_0$ is conjugate to some $m_j^{\pm 1}$, then $m_0$ is a product of
$m_j^{\pm 1}$ with some element of $[\left<m_j\right>,F_m]$, so the following
Lemma 2.9 implies that the Magnus expansion $M(m_0)$ contains at least one
occurence of $x_j$ in every term of positive degree.
Now apply the same lemma inductively to the elements of
$[F_m,_i\left<m_0\right>]$ to obtain that for any
$c\in [F_m,_{k+2}\left<m_0\right>]$, every term of $M(c)-1$ contains $x_j$ with
total degree $\ge k+2$.

On the other hand, it is not hard to verify by induction that any commutator
$c_0$ of weight $r$ in $m_1,\dots,m_m$ has its Magnus expansion $M(c_0)$ of the
form $1+\rho+\xi$ where $\rho=\rho(c_0)$ is a homogenious polynomial of degree
$r$ and every term of $\xi$ has degree $>r$, moreover each $x_i$ appears in
every term of $\rho$ with total degree equal to the number of occurences of
$m_i$ in $c_0$.
Now by the P. Hall Basis Theorem \cite{Hall}, each $c\in\gamma_r F_m$ can be
expressed modulo $\gamma_{r+1} F_m$ as a monomial
$c=c_{i_1}^{e_{i_1}}\dots c_{i_s}^{e_{i_s}}$ in basic commutators of weight
$r$.

Assume on the contrary that $c\in [F_m,_{k+2}\left<m_0\right>]$ but some $c_i$
in the above experession, with $e_i\ne 0$, contains less than $k+2$ occurences
of $m_j$.
We have $\rho(c)=e_{i_1}\rho(c_{i_1})+\dots+e_{i_s}\rho(c_{i_s})$, and since
$\rho(c_{i_1}),\dots,\rho(c_{i_s})$ are linearly independent \cite{Hall},
$\rho(c)$ contains nonzero terms where the total degree of $x_j$ is less than
$k+2$.
On the other hand, since $c\in [F_m,_{k+2}\left<m_0\right>]$, we have shown
above that every term of $\rho(c)$ contains $x_j$ with total degree $\ge k+2$.
This contradiction proves the assertion. \qed
\enddemo

\proclaim{Lemma 2.9} Let $F_m=\left<m_1,\dots,m_m\mid\,\right>$ be the free
group, $P$ be the ring of formal power series with integer coefficients in
noncommuting variables $x_1\dots,x_m$, and $M\:F_m\to P$ be the Magnus
expansion defined by $m_i\mapsto 1+x_i$ (so that
$m_i^{-1}\mapsto 1-x_i+x_i^2-x_i^3+\dots$).
Let $a,b\in F_m$ and $i\in\{1,\dots,m\}$.
If each term of $M(a)-1$ has total degree $\ge d_1$ in $x_i$ and each term of
$M(b)-1$ has total degree $\ge d_2$ in $x_i$, then

a) each term of $M(a^{-1})-1$ has total degree $\ge d_1$ in $x_i$,

b) each term of $M(ab)-1$ has total degree $\ge\min(d_1,d_2)$ in $x_i$,

c) each term of $M([a,b])-1$ has total degree $\ge d_1+d_2$ in $x_i$.
\endproclaim

\demo{Proof of (c)} Denote $l_x=M(x)-1$, then $M([a,b])$ can be rewritten,
using the identities $(l_{a^{-1}}+1)(l_a+1)=1$ and
$(l_{b^{-1}}+1)(l_b+1)=1$, as
$$\multline \hskip -5pt
(l_{a^{-1}}+1)(l_{b^{-1}}+1)(l_a+1)(l_b+1)=
(l_{a^{-1}}+1)l_{b^{-1}}(l_a+1)l_b+(l_{a^{-1}}+1)l_{b^{-1}}(l_a+1)+l_b+1\\
=l_{a^{-1}}l_{b^{-1}}l_al_b+l_{b^{-1}}l_al_b+l_{a^{-1}}l_{b^{-1}}l_b+
l_{b^{-1}}l_b+l_{a^{-1}}l_{b^{-1}}l_a+l_{b^{-1}}l_a+l_{a^{-1}}l_{b^{-1}}+
l_{b^{-1}}+l_b+1\\
=l_{a^{-1}}l_{b^{-1}}l_al_b+l_{b^{-1}}l_al_b+l_{a^{-1}}l_{b^{-1}}l_a+
l_{b^{-1}}l_a-l_{a^{-1}}l_b+1. \qed
\endmultline$$
\enddemo

Theorem 2.8 together with the Hall Basis Theorem yield the following

\proclaim{Corollary 2.10} Every element $\alpha$ of the `$k$-quasi-free' group
$F_m/\mu_k$ (isomorphic to $\G_k(\text{$m$-component unlink})$) can be uniquely
represented as $\alpha=C_1^{e_1}\dots C_t^{e_t}$, $e_i\in\Z$, where
$C_1,\dots,C_t$ are the basic commutators in the fixed set of free generators
of $F_m$ with at most $k+1$ entries of each generator.
\endproclaim

For the rest of this section, let $m_i$ denote a fixed meridian to the
$i^{\text{th}}$ component of $L$ and $l_i$ the corresponding longitude, and let
$\frak l_i$ be a word in $m_1,\dots,m_m$ representing by Lemma 2.2 the class of
$l_i$ in $\pi(L)/\gamma_q$ for some fixed $q$.
Note that it may be impossible to use the same word $\frak l_i=\frak l_i(q)$
for each positive integer $q$.

\proclaim{Corollary 2.11} If $q>m(k+1)$, $\G_k(L)$ can be presented as
$$\G_k(L)=\left<m_1,\dots,m_m \mid [m_1,\frak l_1],\dots,[m_m,\frak l_m],\,
c_1,\dots,c_r,\, d_1,\dots,d_s\right>$$
where $c_1,\dots,c_r$ are all basic commutators of weight $\le m(k+1)$ in
$m_i$'s such that at least $k+2$ entries are occupied by the same symbol,
and $d_1,\dots,d_s$ are all left-normed%
\footnote{See comments to Problem 3.15.}
commutators in $m_i$'s of weight $m(k+1)+1$.
\endproclaim

This follows from 2.8, 2.3 and \cite{Mi2; Theorem 4} where it was proved that
$\pi(L)/\gamma_q$ can be presented as
$\left<m_1,\dots,m_m\mid[m_1,\frak l_1],\dots,[m_1,\frak l_1];\gamma_q\right>$.
The case $k=0$ of Corollary 2.11 improves on a result of Levine
\cite{Le; Proposition 2.11}.

\remark{Remark} In general, a presentation for $\G_k(L)$ can {\it not} be
obtained by adding the relators $c_1,\dots,c_r,d_1,\dots,d_s$ to a presentation
of $\pi(L)$, for $\mu_k\pi(L)$ need not coincide with the normal closure $N$ of
$\phi(\mu_k F_m)$ in $\pi(L)$, where $\phi\:F_m\to\pi(L)$ sends the generators
to the fixed meridians $m_i$.
The latter is obvious e.g\. when $L$ is a split link with knotted components:
the forgetful homomorphism into the fundamental group of a component sends
$\mu_k\pi(L)$ onto the commutator subgroup (since this is the case for
the unlink, which is $k$-quasi-isotopic to $L$), while any $[m_i,m_j]$, hence
$N$, clearly goes to $1$.
\endremark

\bigskip
We recall that the invariant $\bar\mu(i_1\dots i_si)$ \cite{Mi2}, where $s<q$,
is the coefficient $\mu(i_1\dots i_si)$ at $x_{i_1}\dots x_{i_s}$ in the Magnus
expansion (see 2.9) of $\frak l_i$, reduced modulo the greatest common divisor
$\Delta(i_1\dots i_si)$ of all $\bar\mu(j_1\dots j_t)$ where $j_1,\dots,j_t$
runs over all ordered subsequences of $i_1,\dots,i_si$ of length $2\le t\le s$.

Note that in view of cyclic symmetry \cite{Mi2; formula (21)} we may always
assume that if some index occurs `too many' times in a $\bar\mu$-invariant,
one of its entries is in the rightmost place.

\proclaim{Theorem 2.12} Suppose $k\ge 1$.

a) If for some $z$ the sequence $i_1,\dots\hat i_z\dots,i_s$ contains at most
$k+1$ entries of each index, $\bar\mu(i_1,\dots,i_s)$ can be extracted from
$(\G_k,\P_k)$; in particular, it is invariant under $k$-quasi-isotopy.

b) An $m$-component link $L$ has the same $(\G_k,\P_k)$ as the $m$-component
unlink iff all $\bar\mu$-invariants with $\le k+2$ entries of one of the
indices and $\le k+1$ entries of each of the remaining indices vanish.

c) A $\bar\mu$-invariant with at least one index occupying more than $k+2$
entries or at least two indices occupying more than $k+1$ entries each cannot
be extracted from $(\G_k,\P_k)$, unless identically trivial.
\endproclaim

\demo{Proof. (a)} By Theorem 2.8 and Proposition 2.13(i) below, the value of
the $\bar\mu$-invariant on any link can be determined from $\G_k$ together
with the images $(m_i,l_i)$ of meridians and longitudes.
It remains to verify that this value is unchanged under multiplication
of $l_i$ on the left by an element of $N_k(m_i)$.
But this follows from Lemma 2.9 (similarly to \cite{Mi2; proof of (25)}) using
that $2k+1>k+1$, which accounts for the restriction $k>0$. \qed
\enddemo

\demo{(b)} The `only if' part follows from (a).
The `if' part follows from Proposition 2.5 and Proposition 2.13(ii) below. \qed
\enddemo

\demo{(c)} Let $\bar\mu(i_1\dots i_n)$ be the invariant in question.
By \cite{Orr} and Proposition 2.13(iii) below it is possible to construct a
link $L_n$ which has a nontrivial $\bar\mu(i_1\dots i_n)$, trivial
$\bar\mu$-invariants of length $<n$, and $\bar\mu(j_1\dots j_n)=0$ whenever
$(j_1,\dots,j_n)$ is not a permutation of $(i_1,\dots,i_n)$.
By Proposition 2.13(iv), since cyclic symmetry holds for $\bar\mu$-invariants
of arbitrary length \cite{Mi2}, we see that the $\bar\mu$-invariants of $L_n$
of length $n+1$ satisfy all the relations that the $\bar\mu$-invariants of
length $n+1$ of a link with vanishing $\bar\mu$-invariants of length $<n+1$
have to satisfy.
Therefore by \cite{Orr} there exists a link $L_{n+1}'$ with vanishing
$\bar\mu$-invariants of length $<n+1$ and exactly the same $\bar\mu$-invariants
of length $n+1$ as those of $L_n$, but with the opposite sign.
Then by \cite{Kr} any connected sum $L_{n+1}=L_n+L_{n+1}'$ has the same
$\bar\mu$-invariants of length $\le n$ as $L_n$, moreover
$$\bar\mu_{L_{n+1}}(j_1,\dots,j_{n+1})=0\qquad\text{(mod
$\Delta_{L_n}(j_1,\dots,j_{n+1})$)}$$
for any multi-index $(j_1,\dots,j_{n+1})$.

Iterating this process, we obtain a link $L_{m(k+1)+1}$ which has
$\bar\mu(i_1,\dots,i_n)\ne 0$, and all $\bar\mu$-invariants with $\le k+2$
entries of one of the indices and $\le k+1$ entries of each of the remaining
indices trivial.
The assertion now follows from the `if' part of (b). \qed
\enddemo

\proclaim{Proposition 2.13} Let $\Cal I=(i_1,\dots,i_n)$, where $n+1<q$,
denote an unordered multi-index (with possible repeats of entries),
$1\le i_j\le m$, let $\Cal I^\sigma$ denote $\Cal I$ with a fixed order, and
let $\Cal I^c$ denote a basic commutator in the alphabet $m_1,\dots,m_m$ whose
entries are $m_{i_1},\dots,m_{i_n}$ (in some order).
Let $i$ be a fixed index, $1\le i\le m$, and let $E(c_1;w),\dots,E(c_r;w)$
denote the powers of the basic commutators $c_1,\dots,c_r$ of weight $<q$ in
the decomposition of a word $w$ in the free nilpotent group $F_m/\gamma_q$.

(i) The invariants $\bar\mu(\Cal I^\sigma i)$ for all orderings $\sigma$ of
$\Cal I$ can be expressed in the commutator numbers%
\footnote{Beware that these integers are not, in general, invariant even
under ambient isotopy of $L$; see \cite{Le} for a discussion of their
indeterminacies in the case where all indices are distinct.}
$E(\Cal I^c;\frak l_i)$ where $c$ runs over all bracketings of $\Cal I$.

(ii) Conversely, if $\Delta(\Cal I^\sigma i)=0$ for all orderings $\sigma$ of
$\Cal I$, the commutator numbers $E(\Cal I^c;\frak l_i)$ for all bracketings
$c$ of $\Cal I$ can be expressed in the invariants $\bar\mu(\Cal I^\sigma i)$.

(iii) Furthermore, if all $\bar\mu$-invariants of length $\le n$ vanish, then
a complete set of relations between the integer invariants
$E(\Cal J^c;\frak l_j)$, where $\Cal J^c$ runs over all basic commutators
of weight $n$ and $j$ runs over $1,\dots,m$, is given by
$$\sum_{(\Cal J,j)=\Cal K}\sum_c
E(\Cal K^d;[m_j,\Cal J^c])\,E(\Cal J^c;\frak l_j)=0\tag{$*$}$$
for each basic commutator $\Cal K^d$ of weight $n+1$.

(iv) Moreover, relations ($*$) are equivalent to cyclic symmetry of
$\bar\mu$-invariants of length $n+1$.
\endproclaim

(In the last assertion it is understood that the commutator numbers
$E(\Cal J^c;\frak l_j)$ and the $\bar\mu$-invariants are calculated from
arbitrary fixed words $\frak l_i=\frak l_i(m_1,\dots,m_m)$ rather than those
representing the longitudes, and it is assumed that among these formal
$\bar\mu$-invariants, all of length $\le n$ vanish.
Of course, there are other symmetries of $\bar\mu$-invariants, but we see that,
at least in the case where all $\bar\mu$-invariants of length $\le n$ vanish,
all relations between $\bar\mu$-invariants of length $n+1$, except for cyclic
symmetry, come from nontriviality of the cokernel of the Magnus expansion
rather than from geometry.)

\demo{Proof} In the ring $P$ of Lemma 2.9, let $D_i$ denote the subset, which
is easily seen to be a two-sided ideal, consisting of all
$\sum\nu^{h_1\dots h_s}x_{h_1}\dots x_{h_s}$ such that either $s\ge q$ or $s<q$
and $\nu^{h_1\dots h_s}\equiv 0$ (mod $\Delta(h_1,\dots,h_si)$).

It is not hard to see that for every $c\in F_m$, every monomial of
$M(c)-\rho(c)-1$ can be obtained from some monomial of $\rho(c)$, where
$\rho(c)$ denotes the leading part of $M(c)$ (that is, $M(c)=1+\rho(c)+$ terms
of higher degree, where $\rho(c)$ is homogenious) by insertion of at least one
letter and multiplication by a constant.
It follows that $$M(\frak l_i)=M(c_1^{e_1}\dots c_r^{e_r})\equiv
e_1\rho(c_1)+\dots+e_r\rho(c_r)\qquad\text{(mod $D_i$)}.$$
Since the total degree of $x_i$ in every term of $\rho(c_j)$ equals the number
of entries of $m_i$ in $c_j$, and the residue classes
$\bar\mu(\Cal I^\sigma i)$ are well-defined up to congruence of $M(\frak l_i)$
modulo $D_i$, this proves the first assertion.

By the above, we have $\sum_c E(\Cal I^c;\frak l_i)\rho(\Cal I^c)=
\sum_\sigma\mu(\Cal I^\sigma i)x_{\Cal I^\sigma}$ (mod $D_i$, hence
absolutely), and since the leading parts of the Magnus expansions of the
basic commutators of weight $n$ are linearly independent \cite{Hall}, this
proves the second assertion.

It is well-known \cite{Rol} that one of the relations in the Wirtinger
presentation of a link is redundant.
It follows from the precise form of this redundancy (see \cite{Mi2; proof of
Lemma 5} or use the above argument with the tangle) that a product of certain
conjugates of the relators $[m_i,\frak l_i]$ of Milnor's presentation for
$\pi(L)/\gamma_q$ equals $1$ in $F_m/\gamma_q$.
In our situation, where all $\bar\mu$-invariants of length $\le n$ are trivial,
and so all longitudes lie in $\gamma_n$, this implies
$$[m_1,\frak l_1][m_2,\frak l_2]\dots [m_m,\frak l_m]=1\qquad
\text{(mod $\gamma_{n+2}F_m$)}.$$
This identity means that the words $\frak l_1,\dots,\frak l_m$ cannot be chosen
independently of each other, and it follows that their commutator numbers
satisfy relations $(*)$.

From the definition $\gamma_{n+1}F_m=[F_m,\gamma_n F_m]$ it follows that each
of relations $(*)$ contains at least one nonzero coefficient
$E(\Cal K^d;[\Cal J^c,m_j])$.
Therefore we have $N_{n+1}$ of nontrivial linear relations imposed on $mN_n$ of
integer indeterminates, where $N_n$ denotes the number of basic commutators
of weight $n$ (in the alphabet $m_1,\dots,m_m$).
But it was proved in \cite{Orr} that if all $\bar\mu$-invariants of length
$\le n$ vanish, then there are exactly $mN_n-N_{n+1}$ linearly independent
$\bar\mu$-invariants of length $n+1$.
Thus relations $(*)$ are all that hold.

To prove the last assertion, we note that it was shown in
\cite{Mi2; proof of Lemma ~5} that cyclic symmetry of formal
$\bar\mu$-invariants of length $n+1$ is equivalent to the relation
$\rho([m_1,\frak l_1]\dots [m_m,\frak l_m])=0$ (mod $D$),
where $D$ is the ideal formed by all polynomials
$\sum\nu^{h_1\dots h_s}x_{h_1}\dots x_{h_s}$ such that either $s\ge q$ or $s<q$
and $\nu^{h_1\dots h_s}\equiv 0$ (mod formal $\Delta(h_1,\dots,h_s)$).
Now if all formal $\bar\mu$-invariants of length $\le n$ vanish,
we have $\rho([m_1,\frak l_1]\dots [m_m,\frak l_m])=0$ absolutely,
which is equivalent to $[m_1,\frak l_1]\dots [m_m,\frak l_m]=1$
(mod $\gamma_{n+2}F_m$). \qed
\enddemo

\head 3. Milnor's links and bounded Engel elements \endhead

Let $\M\:S^1_1\sqcup S^1_2\emb S^3$ denote the $(k+1)^{\text{th}}$ Milnor
link \cite{Mi2; Fig\. 1} (see Fig. ~5, where $k=3$) where $S^1_1$ is standardly
embedded, and let $T\cong S^1\x D^2$ be the exterior of a regular neighborhood
of $\M(S^1_1)$ in $S^3\but\M(S^1_2)$.
Let us call the embedding $\M|_{S^1_2}\:S^1\to T$ {\it the $k^{\text{th}}$
Milnor curve} in the solid torus $T$ (so that the zeroth is the well-known
Whitehead curve).

\def\epsfsize#1#2{.5\hsize}
\fig 5

We say that an $m$-component link $L'$ is obtained from $L$ by an
{\it elementary $k$-move} if there is a solid torus $T\i S^3\but
L(mS^1)$ with the $k^{\text{th}}$ Milnor curve $C\i T$, so that
$L'$ is the interior connected sum of $L$ and $C$ along some band
$b\:I\x (I,\partial I)\emb (S^3, L(mS^1)\cup C)$ joining the arc
$b(I\x 0)\i C$ with the arc $b(I\x 1)\i L(mS^1)$ and meeting
$\partial T$ in an arc, say, $b(I\x\frac12)$. More precisely,
$$L'(mS^1)=(L(mS^1)\cup C)\but h(I\x\partial I)\cup h(\partial I\x
I).$$ Then {\it fine $k$-quasi-isotopy} is the equivalence
relation generated by elementary $k$-moves and ambient isotopy. It
is straightforward that fine $k$-quasi-isotopy implies
$k$-quasi-isotopy, and fine $0$-quasi-isotopy coincides with
$0$-quasi-isotopy (i.e\. link homotopy).

\proclaim{Proposition 3.1}
If two links are fine $k$-quasi-isotopic, then one can be deformed into another
by Habiro moves \cite{Hab}, corresponding to the graph with $k+2$ univalent
vertices, shown on Fig\. 6.
\endproclaim

\def\epsfsize#1#2{.3\hsize}
\fig 6

\demo{Proof} It suffices to consider links that differ by an
elementary $k$-move. In this case our claim is clear from Fig\. 7,
which demonstrates unlinking of the $(k+1)^{\text{th}}$ Milnor
link by means of Habiro moves, corresponding to the unitrivalent
graphs above. \qed
\enddemo

\def\epsfsize#1#2{\hsize}
\fig 7

This result should be considered in the context of statements on Engel elements
in groups, which we are to discuss next.
The informal idea here is that Habiro graphs from Fig\. 6 pictorially resemble
the Engel structure.
We expect that this observation can be extended to relate other types of Habiro
graphs and various subgroups in link groups.

\proclaim{Theorem 3.2} Let $\delta_k$ denote the subgroup of $\pi(L)$ generated
by $$\{\underbrace{[[\dots[}_{k+2} g,m]\dots,m],m]\mid m\in M, g\in\pi(L)\}.$$
Then $\delta_k$ is normal and coincides with the subgroup of $\pi(L)$ generated
by $$\{[m,\underbrace{m^{.^{.^{.^{m^g}}}}\hskip -8pt}_{k+1}
\hskip 5pt]\mid m\in M, g\in\pi(L)\},$$
and the quotient group $\G_k^+(L)=\pi(L)/\delta_k$ is the finest quotient of
$\pi(L)$, functorially invariant under fine $k$-quasi-isotopy.
\endproclaim

\demo{Proof} First let us notice that $[m,m^g]=[m,m[m,g]]=[m,[m,g]]$, hence
$$[m,\underbrace{m^{m^{\cdot^{\cdot^{\cdot^{m^g}}}}}\!\!\!\!\!}_
{k+1\text{ of }\,m\text{'s}}]=[m,[m,\dots[m,g\underbrace{]\dots]]}_{k+2}.$$
$$\multline\text{Also}\qquad
[m,[m,\dots[m,g]\dots]]=[[\dots[g,m]^{-1}\dots,m]^{-1},m]^{-1}\\ =([([\dots
([g^{m^{k+2}},m]^{-1})^{m^{-1}}\dots,m]^{-1})^{m^{-1}},m]^{-1})^{m^{-1}}\\
=[[\dots[g^{m^{k+2}},m^{-1}]\dots,m^{-1}],m^{-1}].
\endmultline$$
This means, using the notation $m^g_0=g$, $m^g_{k+1}=m^{m^g_k}$, that
$\delta_k$ is generated by $$\{[m,m^g_{k+1}]\mid m\in M, g\in\pi(L)\}.$$
A standard computation based on \cite{MR1; proof of Theorem 3.2} shows that
$\G^+_k$ is functorially invariant under fine $k$-quasi-isotopy.

So (cf\. \cite{MR1; proof of Theorem 3.7}) it remains to show
that for any link $L_0$, any meridian $m$ of $L_0$, and any $g\in\pi(L_0)$
there exists a fine $k$-quasi-isotopy $L_t$ with a single self-intersection of
a component such that the link $L_{s-\eps}$ immediately preceding the singular
link $L_s$, has meridians $m$ and $m^g_{k+1}$ sharing the same stem and winding
around the arcs $L_{s-\eps}(l_1)$ and $L_{s-\eps}(l_2)$ respectively, where
$l_1\sqcup l_2\i S^1_j$ is a small neighborhood of the preimage of the double
point of $L_s$.

To see this, first in the special case where $L_0$ is the two-component unlink
with image in the plane $\R^2\x\{0\}\i\R^3$ and $m,g$ are the Wirtinger
generators of $\pi(L_0)$, let us consider the evident straight line homotopy
with a single self-intersection from the Milnor link $\M_{k+1}$ to the link
$M_{k+1}$, which is isotopic to the unlink and differs from $\M_{k+1}$ only in
that the clasp is `open'.
One verifies that the meridians $m$ and $m^g_{k+1}$ of the link $M_{k+1}$ are
exactly as required, and our special case follows.

The general case follows by applying the special case locally.
We leave details to the reader. \qed
\enddemo

In fact, the difference between the subgroups $\mu_k$ and $\delta_k$ is similar
to the difference between subgroups $\nu_{k+2}$ and $\eps_{k+2}$ defined below,
which are of some importance in combinatorial group theory.

\subhead Engel groups and $\eps_n$ \endsubhead
A group $G$ is called an {\it Engel group} if
$$\underbrace{[[\dots[}_ng,h]\dots,h],h]=1$$ for each $g,h\in G$ and
some $n=n(g,h)$, and an {\it $n$-Engel group} if moreover $n$ can be chosen
independent of $g$ and $h$ \cite{Gr1}, \cite{Ku}, \cite{KMe}, \cite{Rob}.
An element $h$ of an arbitrary group $G$ is called a {\it (left) Engel element}
if the equality above holds for each $g\in G$ and some $n=n(g,h)$, and a {\it
(left) $n$-Engel element} if in addition $n$ can be chosen the same for all $g$
\cite{Gr2}, \cite{Ku}, \cite{Rob}.
For any group $G$ we denote by $\eps_n G$, $n\ge 1$, its $n^{\text{th}}$
{\it Engel subgroup}, generated by the set
$$\{\smash{\underbrace{[[\dots[}_n}g,h]\dots,h],h]\mid g,h\in G\}.$$
\bigskip\smallskip\noindent
Since $[g,h]^f=[g^f,h^f]$, this subgroup is normal.
It is easy to see that a group $G$ is $n$-Engel iff $\eps_n G=1$.
Clearly, $\eps_1 G=\gamma_2 G$ for any $G$.
Notice that $\eps_2 G$ is generated by $\{[h,h^g]\mid g,h\in G\}$ since
$\eps_2 G$ is normal and
$$[[g,h],h]=[[h,g]^{-1},h]=[h,[h,g]]^{[g,h]}=[h,h[h,g]]^{[g,h]}=
[h,h^g]^{[g,h]}.$$
In this form $\eps_2$ looks more similar to $\gamma_2$ rather than $\gamma_3$.

\remark{Remark} The above observation that $[g,h,h]$ is conjugate with
$[h,h^g]$ can be pushed one step further, that is, in any group, $[g,h,h,h]$
is conjugate with $[h,h^{h^g}]^{-1}$.
Indeed, notice that $[h,h^g]=[h,h[h,g]]=[h,[h,g]]$, hence
$[h,h^{h^g}]=[h,[h,[h,g]]]$.
So $$\multline
[g,h,h,h]=[[h,[h,g]]^{[g,h]},h]=[[h,[h,g]],h^{[h,g]}]^{[g,h]}
=[[h,[h,g]],h[h,[h,g]]]^{[g,h]}\\
=[[h,[h,g]],h]^{[h,[h,g]][g,h]}=[h,[h,[h,g]]]^{-[g,h]^h}.
\endmultline$$
\endremark

\example{Example 3.3}
Here are two simple cases where $\eps_2 G\ne\gamma_2 G$.
The group
$$Q_8=\left<i,j,k\mid i^2=j^2=k^2=ijk=(ijk)^{-1}\right>$$ is the union of its
subgroups $\left<i\right>$, $\left<j\right>$, $\left<k\right>$ so that any two
elements, conjugate in $Q_8$, lie in one of them, which means that
$\eps_2 Q_8=1$, whereas $\gamma_2 Q_8=\left<ijk\right>$.
Alternatively, let $D_4$ be the $2$-Sylow subgroup of the symmetric
group $S_4$, generated by the permutations $(12)$ and $(13)(24)$.
Then $D_4$ is the union of its subgroups
$$\gather H_1=\{(12),(34),(12)(34),1\},\\
H_2=\{(13)(24),(14)(23),(12)(34),1\},\\
H_3=\{(1324),(1423),(12)(34),1\}.\endgather$$
Any two permutations, conjugate in $D_4$, should lie in one of these Abelian
subgroups.
On the other hand, $[D_4,D_4]=[H_i,H_j]=\{(12)(34)\}$, $i\neq j$.
\endexample

\example{Example 3.4} Let us see that $\eps_2 F_m\ne\gamma_3 F_m$ for $m\ge 3$,
where $F_m=\left<a_1,\dots,a_m\right>$ denotes the free group.
It suffices to prove that $\eps_2 F_3/\gamma_4 F_3\ne\gamma_3 F_3/\gamma_4F_3$.
The basic commutators of weight $3$ in $F_3$ are
$$[[b,a],a]\quad [[b,a],b]\quad [[b,a],c]\quad [[c,a],a]\quad [[c,a],b]\quad
[[c,a],c]\quad [[c,b],b]\quad [[c,b],c]].$$
Let us consider the quotient $Q$ of $\gamma_3 F_3/\gamma_4 F_3$ by the
subgroup $$(\gamma_3\left<a,b\right>/\gamma_4 F_3)\oplus
(\gamma_3\left<a,c\right>/\gamma_4 F_3)\oplus
(\gamma_3\left<b,c\right>/\gamma_4 F_3).$$
Now $[[h,g],g]$ can be written modulo $\gamma_4 F_3$ as
$[[a^ib^jc^k,a^lb^mc^n],a^lb^mc^n]$ for some $i,j,k,l,m,n\in\Z$, and
$$\multline [[a^ib^jc^k,a^lb^mc^n],a^lb^mc^n]\underset Q\to=
[[a,b^mc^n],b^mc^n]^i[[b,a^lc^n],a^lc^n]^j[[c,a^lb^m],a^lb^m]^k\underset Q\to=\\
[[a,b],c]^{inm}[[a,c],b]^{inm}[[b,a],c]^{jln}[[b,c],a]^{jln}[[c,a],b]^{klm}
[[c,b],a]^{klm}\underset\gamma_4\to=\\ [[b,a],c]^{-inm}[[c,a],b]^{-inm}
[[b,a],c]^{2jln}[[c,a],b]^{-jln}[[b,a],c]^{-klm}[[c,a],b]^{2klm},\endmultline$$
where the last line uses the Hall--Witt identity
$[[a,b^{-1}],c]^b[[b,c^{-1}],a]^c[[c,a^{-1}],b]^a=1$ in the form
$[[a,b],c][[b,c],a][[c,a],b]\underset\gamma_4\to=1$.

Since $Q$ is generated by two basic commutators, the quotient $Q'$ of $Q$ by
the subgroup $\left<[[b,a],c][[c,a],b]\right>$ is infinite cyclic
generated by $[[c,a],b]\underset Q'\to=[[b,a],c]^{-1}$.
From the above we see that the image of $\eps_2 F_3/\gamma_4 F_3$ in $Q'$
is the subgroup $3Q'\ne Q'$, which proves the assertion.
\endexample

\remark{Remarks. (i)}
Notice that $\eps_2 F_2=\gamma_3 F_2$ since $\eps_2 F_2$ contains
$\gamma_4 F_2$ \cite{Rob}, meanwhile both basic commutators of
weight $3$ in $F_2$ lie in $\eps_2 F_2$, namely:
$$\gather [[b,a],a]\in \eps_2 F_2\\
[[b,a],b]=[[a,b]^{-1},b]=([[a,b],b]^{-1})^{[b,a]}\in\eps_2 F_2.\endgather$$
\endremark

\remark{(ii)} It is not hard to show, using the above methods, that
$\eps_3 F_2/\gamma_5 F_2$ is a subgroup of index two in
$\gamma_4 F_2/\gamma_5 F_2$, and therefore $\eps_3 F_2\ne\gamma_4 F_2$.
\endremark

\remark{(iii)} From the above argument we also see that the product of two
$2$-Engel elements is not necessarily a $2$-Engel element.
\endremark
\bigskip

Notice that to prove that a group, generated by $(k+2)$-Engel elements
(in particular, $\G_k^+(L)$ for any link $L$) may be non-nilpotent, one has
to deal with at least one of the following two group-theoretic problems,
dating back to 1950s \cite{Gr1}, \cite{Ba2}, \cite{Pl1}, which to the best
of the authors' knowledge are still unsolved (cf\. \cite{KMe}, \cite{Rob},
\cite{BM}, \cite{Va}, \cite{Ko}):

(i) whether every finitely generated group where each element is bounded Engel,
must be nilpotent, and

(ii) whether the set of bounded Engel elements of every group is a subgroup.

\noindent
These are modifications of the two famous problems: whether finitely generated
$n$-Engel groups are nilpotent, and whether the set of Engel elements is
always a subgroup.
There are simple examples of non-nilpotent infinitely generated $n$-Engel
groups \cite{Rob; Ex\. 12.3.1}, \cite{KMe} and an example, due to Golod (1964),
of a non-nilpotent finitely generated Engel group \cite{KMe}.

\proclaim{Theorem 3.5}
a) Consider $\G^{++}_k(L)=\pi(L)/\eps_{k+2}$, the `homogenious' quotient of
$\G^+_k(L)$.
Then $\G^{++}_1$ is nilpotent, all epimorphic images of $\G^{++}_2$ of exponent
$5$ are nilpotent, and for any $k$, all residually finite epimorphic images of
$\G^{++}_k$ are nilpotent.

b) Any epimorphic image of $\G^+_k(L)$ with all Abelian subgroups finitely
generated or with all normal closures of single elements finitely generated,
is Engel.

c) Any solvable or Noetherian (in particular, any finite) epimorphic image of
$\G^+_k(L)$ is nilpotent.
\endproclaim

Recall that a group is {\it Noetherian} if all its subgroups are finitely
generated.
The third part of (a) can be compared with the fact \cite{Hem} that
the fundamental group of any link is {\it residually finite} (i.e\. admits,
for each non-trivial element, an epimorphism onto a finite group which does
not trivialize this element).

\demo{Proof. (a)} By the definition, $\G^{++}_k$ is a $(k+2)$-Engel group.
Now finitely generated $3$-Engel groups are known to be nilpotent \cite{Hei}.
The second and third assertions follow from the results of Vaughan-Lee
\cite{Va} and Wilson \cite{Wi} that all $4$-Engel groups of exponent $5$ and
all finitely generated residually finite $l$-Engel groups are nilpotent. \qed
\enddemo

\remark{Remark} For finite, rather then residually finite, epimorphic image
it suffices to use the easier result that any finite Engel group is nilpotent
(see \cite{Rob} and references in \cite{Gr1} and \cite{Ku}).
\endremark

\demo{(c)} The image $\bar m_i$ of each meridian $m_i$ in the quotient
$\G^+_k(L)$ is a $(k+2)$-Engel element.
This holds if we proceed further to the given solvable (or Noetherian) quotient
$Q$ of $\G^+_k(L)$, that is, $\bar{\bar m}_i\in Q$ is a $(k+2)$-Engel element.
Since $Q$ is solvable (resp\. Noetherian), by \cite{Gr2} (resp\. by \cite{Ba2})
all its Engel elements form a normal subgroup.
But the $m_i$'s normally generate $\pi(L)$, hence the $\bar{\bar m}_i$'s
normally generate the whole $S$, thus $Q$ is an Engel group.
Now every Noetherian (resp\. finitely generated solvable) Engel group is
nilpotent by \cite{Ba2} (resp\. by \cite{Gr1}, see also \cite{Rob}) which
completes the proof. \qed
\enddemo

\demo{(b)} This is analogous to the above argument, using the result of Plotkin
\cite{Pl2} that in a group with all Abelian subgroups (or all normal closures
of elements) finitely generated the set of Engel elements is a normal subgroup.
\qed
\enddemo

\subhead Baer groups and $\nu_n$ \endsubhead
Let $\nu_n G$, $n\ge 1$, denote the subgroup generated by
$$\bigcup_{g\in G}\underbrace{[[\dots[}_n G,\left<g\right>]\dots,\left<g\right>
],\left<g\right>],$$
then $\nu_n G$ is normal in $G$.
The group $G$ is called a {\it Baer group} \cite{Ba1}, \cite{Rob}, if for each
$g\in G$, the cyclic subgroup $H=\left<g\right>$ is subnormal%
\footnote{A subgroup $H$ of a group $G$ is said to be {\it subnormal} in $G$
if there exists a finite chain of subgroups $H=H_0\i H_1\i\dots\i H_d=G$ such
that each $H_i$ is normal in $H_{i+1}$.}
in $G$.
By Lemma 2.6, $\nu_n G$ coincides with $\left<\bigcup_{g\in G}[\left<g\right>,
\left<g\right>_{n-1}^G]\right>$, so $\nu_n G=1$ implies
$\left<g\right>^G_n=\left<g\right>$ (but not vice versa, as seen,
say, for $n=1$ and the group $D_4$).
In particular, if $\nu_n G=1$ for some finite $n$, then $G$ is a Baer group.

\proclaim{Theorem 3.6} $\nu_n G=1$ iff $\left<g\right>^G$ is nilpotent of class
$n-1$ for each $g\in G$.
\endproclaim

\demo{Proof} The `if' part holds by the definition, since
$[G,\left<g\right>]\i\left<g\right>^G$.
To prove the converse, we may assume that $G$ is finitely generated.
(Indeed, we need to show that for any $g\in G$, any commutator
$c=[h_1,\dots,h_{n-1}]$, where each $h_i$ is a product of conjugates of $g$
by some $h_{ij}\in G$ and their inverses, lies in $\nu_n G$.
But $c\in\gamma_{n-1}\left<g\right>^H$ and $\nu_n H\i\nu_n G$, where $H$ is
the subgroup generated by $g$ and all $g_{ij}$'s.)
Now the proof is analogous to the proof of the inclusion `$\supset$' in
Theorem 2.7, taking into account that finitely generated Baer groups are
nilpotent \cite{Ba1} or, alternatively, that
$\nu_n F_m\supset\mu_{n-2} F_m\supset \gamma_{m(n-1)+1}$ by Theorem 2.1. \qed
\enddemo

Evidently, $\eps_1 G=\nu_1 G=\gamma_2 G$ and in general
$\eps_n G\i\nu_n G\i\gamma_{n+1} G$.
However, it is not obvious that $\nu_n G$ may differ from $\eps_n G$ for some
$n$ and $G$.
The possible difference between $\eps_2 G$ and $\nu_2 G$ could be that
the latter contains all the elements of the form
$$[[h_1,g^{m_1}]^{\alpha_1}[h_2,g^{m_2}]^{\alpha_2}\cdots
[h_r,g^{m_r}]^{\alpha_r},g^m]$$
where all $\alpha_i=\pm 1$, $m_i,m\in\Z$, $h_i\in G$ and $g\in G$, while the
former does not obviously contain them, unless $r=1$ and $\alpha_1=m_1=m=1$.
However it turns out that $\eps_2 G$ and $\nu_2 G$ happen to be identical.
Indeed, since $\eps_2 G$ is normal, it contains all elements of type
$([[h,g],g]^{\pm 1})^c$, meanwhile $$[[h_1,g^{m_1}]^{\alpha_1}\cdots [h_r,g^{m_r}]^{\alpha_r},g^m]=
\prod_{i=1}^{|m|}\prod_{j=1}^r\prod_{k=1}^{|m_j|}
([[h_j,g],g]^{\text{sgn}(m_j)\alpha_j\text{sgn}(m)})^{c_{ijk}}$$
(for certain $c_{ijk}\in G$) due to the commutator identities
$[ab,c]=[a,c]^b[b,c]$ and $[x^{-1},y]=([x,y]^{-1})^{x^{-1}}$.
Notice that, using the identity $x^y=x[x,y]$, the argument above generalizes to
prove that $\nu_n G/\gamma_{n+2} G=\eps_n G/\gamma_{n+2} G$ (cf\. \cite{Mo}).

Furthermore, it was shown in \cite{KK} that $\eps_3=\nu_3$ in any group
(here and below $\nu_n$ is understood to be the subgroup generated by the
$(n-1)^{\text{th}}$ lower central terms of the normal closures of elements).
Around 1980, Gupta, Heineken and Levin constructed groups with $\eps_4\ne\nu_4$
\cite{GL}.
For $k\ge 5$ there are also examples of $k$-Engel groups with $\nu_n\ne 1$ for
any finite $n$ \cite{GL}, but no such finitely generated groups seem to be
known.
The groups in these examples are solvable, although $\eps_n G=\nu_n G$ for
metabelian $G$ \cite{KMo} (see \cite{Mo} for further results in this
direction).

\proclaim{Theorem 3.7} Any quotient of the metabelianized fundamental group,
functorially invariant under $k$-quasi-isotopy, is functorially invariant
under fine $k$-quasi-isotopy.
\endproclaim

Speaking informally, the Alexander module ``cannot distinguish''
$k$-quasi-isotopy from fine $k$-quasi-isotopy.

\demo{Proof}
For any elements $b_1,\dots,b_n$ of a metabelian group $G$ and any $a\in G'$,
one has $[a,b_1,...,b_n]=[a,b_{\sigma(1)},...,b_{\sigma(n)}]$, where $\sigma$
is any permutation, see e.g\. \cite{Ne}.
As noted in the proof of Theorem 2.7, $\mu_k$ is generated by the left-normed
commutators where one of the fixed meridians $m_1,\dots,m_m$ occurs at least
$k+2$ times as an entry.
The above identity can be used to transform any such commutator, modulo
$\pi(L)''$, to a commutator (of weight $1$ perhaps) where one of the entries
is either
$$x:=[g,m_i,\underbrace{m_i\dots,m_i}_{k+1}]\qquad\text{or}\qquad
y:=[m_i,g,\underbrace{m_i\dots,m_i}_{k+1}]$$
for some $g\in\pi(L)$.
Now notice that since $[a^{-1},b]=[a,b]^{-1}[[a,b]^{-1},a^{-1}]$, for any
element $b$ of a metabelian group $G$ and any $a\in G'$, one has
$[a^{-1},b]=[a,b]^{-1}$.
Repeated use of this identity shows that $x=y^{-1}$ modulo $\pi(L)''$.
Thus the images of $\mu_k$ and $\delta_k$ in the metabelianization coincide.
\qed
\enddemo

\proclaim{Conjecture 3.8}
For $k\ge 2$, but not for $k=1$, difference between $k$-quasi-isotopy and
fine $k$-quasi-isotopy can be detected by the higher (non-commutative)
Alexander modules of Cochran et al.
\endproclaim

\proclaim{Theorem 3.9}
A group $G$, generated by two $3$-Engel elements (in particular, $\G_k^+(L)$
for any $2$-component link $L$ with $\pi(L)$ generated by two meridians), is of
length $\le 5$, i.e\. $\gamma_5 G=\gamma_6 G=\dots$
\endproclaim

\remark{Remark}
A group, generated by two $2$-Engel elements, is nilpotent of class $2$,
since the two basic commutators of weight $3$ in $F_2=\left<a,b\mid\ \right>$
normally generate $\gamma_3 F_2$ (cf\. Lemma 2.3), while they can be expressed
as
$$\align
[[b,a],a]&=[b,a,a];\\
[[b,a],b]&=[[a,b]^{-1},b]=([a,b,b]^{-1})^{[b,a]}.
\endalign$$
\endremark

\demo{Proof}
It suffices to prove that the group
$$G=F_2/\delta_1 F_2=\left<a,b\mid [x,a,a,a]=[x,b,b,b]=1\ \,
\forall x\in G\right>$$ has length $\le 5$.
In other words, the image of $\gamma_5 F_2$ in the quotient of $F_2$ by
$(\gamma_6 F_2)(\delta_1 F_2)$ must be trivial.
The abelian group $\gamma_5 F_2/\gamma_6 F_2$ is freely generated by the basic
commutators
$$[b,a,a,a,a];\ \ [b,a,b,b,b];\ \ [b,a,a,a,b];\ \ [b,a,a,b,b];\ \
[b,a,a,[b,a]];\ \ [b,a,b,[b,a]].$$

It is immediate that the first three commutators are trivial modulo
$\delta_1$.

\proclaim{Claim 3.10} In any group, $[h,g,h,h]$ is conjugate to
$[g,h,h,h]^{-1}$.
\endproclaim

\demo{Proof} Indeed,
$$\multline [h,g,h,h]=[[[g,h]^{-1},h],h]=[[h,[g,h]]^{[h,g]},h]
=[[h,[g,h]],h^{[g,h]}]^{[h,g]}\\
=[[h,[g,h]],h[h,[g,h]]]^{[h,g]}=[[h,[g,h]],h]^{[h,[g,h]][h,g]}
=[[g,h,h]^{-1},h]^{[h,[g,h]][h,g]}\\
=[g,h,h,h]^{-[h,[g,h]]^2[h,g]}. \qed
\endmultline$$
\enddemo

\proclaim{Claim 3.11} In any group, $[g,h,h,[h,g]]$ and
$[h,g,h,[h,g]]$ lie in the normal closure of two commutators $[g,h,h,h]$ and
$[g^{-1},h,h,h]$.
\endproclaim

\demo{Proof} From 3.10, $[[g^{\pm 1},h]^{-1},h,h]$ lies
in the normal closure of $[[g^{\pm 1},h],h,h]$.
Now let us modify:
$$\multline
[[g^{-1},h]^{\pm 1},h,h]=[[g,h]^{\mp g^{-1}},h,h]
=[[g,h]^{\mp 1},h^g,h^g]^{g^{-1}}\\
=[[g,h]^{\mp 1},h[h,g],h[h,g]]^{g^{-1}}
=[[[g,h]^{\mp 1},h]^{[h,g]},h[h,g]]^{g^{-1}}\\
=[[g,h]^{\mp 1},h,[h,g]h]^{[h,g]g^{-1}}
=([[g,h]^{\mp 1},h,h][[g,h]^{\mp 1},h,[h,g]]^h)^{g^{-h}}.
\endmultline$$
Consequently $[[g,h]^{\mp 1},h,[h,g]]=[[g,h]^{\mp 1},h,h]^{-h^{-1}}
[[g^{-1},h]^{\pm 1},h,h]^{h^{-1}g}$ and the statement follows. \qed
\enddemo

\remark{Observation 3.12} The statement of Claim 3.11 also applies to
$[g,h,h,[g,h]]$ and $[h,g,h,[g,h]]$, since they are conjugate to
$[g,h,h,[h,g]]^{-1}$ and $[h,g,h,[h,g]]^{-1}$ respectively.
\endremark

\medskip
From 3.11 and 3.12 we conclude that the basic commutators $[b,a,a,[b,a]]$ and
$[b,a,b,[b,a]]$ also lie in $\delta_1$.
Therefore it remains to prove the following.

\proclaim{Lemma 3.13} $[b,a,a,b,b]\in(\gamma_6 F_2)(\delta_1 F_2)$.
\endproclaim

\demo{Proof} Let us modify:
$$\multline
[a,b,b,b]^a=[a,b^a,b^a,b^a]=[a,b[b,a],b^a,b^a]=[[a,[b,a]][a,b],b^a,b^a]\\
=[[a,[b,a],b^a]^{[a,b]}[a,b,b^a],b^a]
=[[a,[b,a],b^a]^{[a,b]},b^a]^{[a,b,b^a]}[a,b,b^a,b^a].
\endmultline$$
Using Claim 3.11,
$$\multline
[a,b,b^a,b^a]=[[a,b,b[b,a]],b[b,a]]=[[a,b,b]^{[b,a]},b[b,a]]
=[[a,b,b],[b,a]b]^{[b,a]}\\
=([a,b,b,b][a,b,b,[b,a]]^b)^{[b,a]}\in\delta_1.
\endmultline$$
Using Observation 3.12,
$$\multline
[[a,[b,a],b^a]^{[a,b]},b^a]=[[a,[b,a],b[b,a]]^{[a,b]},b[b,a]]
=[[a,[b,a]]^{[a,b]},[b,a]b,b[b,a]]\\
=[[a,b,a],[b,a]b,b[b,a]]=[[a,b,a,b][a,b,a,[b,a]]^b,b[b,a]]\underset\delta_1\to=
[a,b,a,b,b[b,a]]\\
=[a,b,a,b,[b,a]][a,b,a,b,b]^{[b,a]}\underset\gamma_6\to=[a,b,a,b,b].
\endmultline$$
Finally, we see that
$$\multline
[a,b,a,b,b]=[[b,a,a]^{-[a,b]},b,b]=[[b,a,a]^{-1},b^{[b,a]},b^{[b,a]}]^{[a,b]}\\
=[[b,a,a]^{-1},b[b,[b,a]],b[b,[b,a]]]^{[a,b]}\underset\gamma_7\to=
[[b,a,a]^{-1},b,b]=[[b,a,a,b]^{-[a,[b,a]]},b]\\
\underset\gamma_8\to=[b,a,a,b,b]^{-1}. \qed \qed
\endmultline$$
\enddemo
\enddemo

\remark{Remark} Let us sketch an alternative proof of Theorem 3.9 avoiding
lengthy calculations of Lemma 3.13.
By Fitting's Theorem, $G/(\gamma_3 A)(\gamma_3 B)$ is nilpotent, where $A$ is
the normal closure of $a$ and $B$ is the normal closure of $b$.
We want to show that, modulo $\gamma_6$, every element of
$(\gamma_3 A)(\gamma_3 B)$ lies in $\delta_1$.
It is easy to see that $(\gamma_3 A)(\gamma_3 B)$ is normally generated by the
commutators $[m,m^g,m^h]$ where $m\in\{a,b\}$ and $g,h\in G$.
Now $[a,a^g,a^h]$ is conjugate to $[a,g,a,a^h]=[a,g,a,a]^{[a,h]}[a,g,a,[a,h]]$,
and by Claim 3.10 the first factor is in $\delta_1$, while the second factor
can be written, modulo $\gamma_6$, as a power of $[a,b,a,[a,b]]$ which by
Claim 3.11 is also in $\delta_1$.
\endremark

\remark{Remark} An earlier version of this paper contained the conjecture that
the quotient $F_2/\delta_1$ is not nilpotent.
We were informed by A. Abdollahi that he recently disproved this conjecture.
In more detail, he proved that $F_2/\delta_1$ is metabelian and deduced that
it is nilpotent of class 4.
This implication also follows from the above computations and the fact that
$\gamma_5F_2$ is the normal closure of basic commutators of weight $5$
\cite{JGS}, or alternatively from Theorems 3.7 and 2.1.
\endremark

\proclaim{Conjecture 3.14}
a) The group $F_k/\delta_1$ is nilpotent of class $2k$;

b) The group $F_2/\delta_2$ is not nilpotent.
\endproclaim

The fact that $F_2/\delta_1$ is nilpotent implies that $\G_1^+(L)$ is not
a complete invariant of fine 1-quasi-isotopy.
Indeed, $\G_1^+$ was shown to have length $2$ for the links from Fig\. 1,
whence it is has to be isomorphic to $\Z\oplus\Z$.
Moreover, the alternative proof of Theorem 3.9 above shows that
$\G_1(L)=\G_1^+(L)$ for any $2$-component link whose fundamental group is
generated by two meridians.

In conclusion we state a problem, which is of some interest in light of
the computations throughout the paper.

\proclaim{Problem 3.15} Let $\left<X\mid R\right>$ be a finite presentation of
a group $G$, where $R$ is some collection of basic commutators (perhaps of
different weights) in some Hall basis.
Is it true that $G$ is residually nilpotent?
\endproclaim

The affirmative answer implies (by the collecting process
\cite{Hall; (11.1.4)}) that $\gamma_nF_k$ coincides with the normal closure of
the set of basic commutators of weight $n$, which is unsolved for $n>5$ and
for $n=4$, $k>2$ \cite{Ko}, \cite{JGS}.
Note that Problem 3.15 is also related to some problems from \cite{Co}, in
particular, whether the fundamental group of the Whitehead link $\Cal W_1$
is residually nilpotent.
It is not hard to show that the group $G$ is residually nilpotent in
the simplest case $\left<a,b\mid [b,a,a]=1\right>$, but already for the series
of cases $\left<a,b\mid [b,a,a,\dots,a]=1\right>$ the problem is unsolved.

\enddocument
\end